\documentclass[12pt]{article}
\usepackage{amsmath, amsfonts, amssymb, amsthm}
\usepackage{epsf}

\usepackage{amsmath,amsfonts,amssymb,amsthm,graphicx,color} \parskip1pt

\def\S{\d s}
\def\s{\'s}

\def\A{\= A}
\def\a{\= a}

\def\u{\= u}

\def\suls{{\'Sulvas\u tras}}

\def\sul{{\'Sulvas\u tra}}

\def\bsl{{Baudh\a yana}}
\def\asl{{\A pastamba}}
\def\ksl{{K\a ty\a yana}}
\def\msl{{M\a nava}}
\def\n{\. n}
\def\nd{\d n}
\def\md{\d m}
\def\dd{\d d}
\def\sd{\d s}
\def\h{\d h}
\def\D{\d d}
\def\r{\d r}
\def\h{\d h}

\begin{document}  
\begin{center}
{{\Large {\bf Cognition of the circle in ancient 
 India}}}

\bigskip

{S.G. Dani}

\end{center}
\medskip

Next only to the rectilinear figures such as triangles and rectangles,
the circle is the simplest geometrical figure that  would have
touched human life even in the primitive stages, especially
after the advent of the wheel. Apart from the everyday secular aspects
of life, the circle gained significance in ritual and spiritual
respects. Considerable understanding was acquired over a period in the
ancient times, concerning various geometric features
of the circle. Progress in understanding the circle may  be readily correlated with progress of human 
civilization in general. Our overall knowledge of history across ancient 
cultures has
many limitations, in terms of source material and means of interpretation. Nevertheless, in the Indian context we are endowed 
with information from various sources such as \suls, the
Jaina compositions, works from the mathematical astronomy tradition 
starting with \A ryabha\d ta, and 
finally the Kerala school of mathematics, from different periods in history,   
that give an interesting perspective on how ideas developed on the
issue. 

\bigskip
\noindent{\bf {\large The \sul\ period}}

\medskip
The \suls\ are compositions concerned with construction of altars
({\it Vedi}) and fire platforms ({\it citi}) for the {Vedic}
rituals.\footnote{Performance of {\it yajna}s, fire rituals, in pursuance of  
material and/or spiritual benefits is one of the dominant features of the {\it Vedic} 
civilisation. It  involved both  the royalty as well as laity from the priestly class of the time. There are detailed 
prescriptions, about specific  {\it yajna}s to be performed for various  objectives, as well as  the procedures to be followed.} Apart from 
elaborate instruction on laying of bricks (of simple rectilinear shapes, such as squares, triangles etc.) to achieve approximations to various elaborate shapes such as 
falcons, tortoise, wheel, etc., the compositions also include enunciation of various geometric principles, geometric constructions with practical or theoretical import, etc., thus providing us a glimpse of the mathematical knowledge at that time. There were many \suls, of which \bsl, \asl, \msl\ and \ksl\  
\suls\ are especially noted for their 
significance from a mathematical point of view. The  period of the \suls\ is somewhat uncertain, 
as there are no internal clues in the compositions other than their style 
and language, but there now seems to be a general consensus among scholars that they 
belong to the period from about 
800 BCE to 200 BCE, \bsl\ being the oldest and \ksl\ the latest. 
For further general details the reader is referred to
\cite{SB}, \cite{Ku}, \cite{D}, and the references cited there - here we
shall focus on the specific theme at hand, concerning the circle. 

One of the simplest questions that one can think of about the circle 
is the ratio of its circumference to the diameter. As in other ancient 
cultures, in 
ancient India also this ratio was believed to be $3$. 
In the context of the {\it Vedic} tradition this is reflected in an indirect 
reference in the \bsl\ \sul\ in the statement
``The pits for the sacrificial posts are 1 {\it pada} in diameter, 3
{\it pada}s in circumference.'' (\bsl\ \sul\ 4.15, see \cite{SB}); {\it pada}, which literally means
foot, was a 
measure of length equivalent to about $28$ cm. The second part
of the statement is evidently meant as an elaboration/clarification of the first part, but provides us a 
clue that they considered the circumference to be 3 times the diameter. 

The choice of the value 3 for the ratio would today seem quite surprising, as
one would expect that many everyday experiences could have suggested
the value to be a little more.  The following seems to me to be a
plausible explanation for this (which does not seem to have appeared in 
literature before): the idea of the ratio being 3 
dates back to the time when mankind was yet to think in terms of fractions
(except 
perhaps for ``half'', which may have meant a substantial portion of the whole,
rather than its precise value as we understand it) and developed into a belief (perhaps linked with religion). The ratio was assigned the value 
 3 in the sense that it is not, say 2 or 4, or even ``three hand a half''. The belief, once it was
rooted deeply, was not reviewed for a
long time, even after fractions became part of human thought process. 
While our encounter with the circle, especially in the context of wheels, is
over 5000 years old, fractions seem to have appeared on the scene in a
serious way, 
in Indian as well as Egyptian cultures, substantially later, possibly only in the first millennium BCE.   
The difference between the actual ratio and 3 is small enough not to come in to
serious conflict with  
everyday experience to warrant doubting an accepted proposition, which furthermore 
may have the backing of religious authority, and an appeal on account of universal significance 
associated with the number three. Also, while
for a first-time determination of an entity one typically 
avails of prevalent  techniques of any given time, a belief  often 
remains untested until coming in conflict  with another idea or
experience. 

The \msl\  \sul\ however breaks out of the mould, and we encounter
the following: 

\begin{center}
{\tt vi\S kambha\h pa${\tilde {\tt n}}$cabh\a ga\s ca vi\S kambhastrigu\nd a\s
  ca ya\h} $|$ 

{\tt sa ma\nd \D alaparik\S epo na v\a
  lamatiricyate} $||$

\smallskip
\hfill (\msl\ \sul\ 10.3.2.13)
\end{center}

\begin{quote}
A fifth of the diameter and three times
  the diameter, that is the circumference of the circle, not even a
  hair-breadth remains.
\end{quote}

Apparently over the years it was recognised that the ratio is indeed a little 
larger than $3$. \msl\  seems to have taken a leaf out of this and came up
with a better estimate. The exultation over it is striking!

\medskip
Unlike the  circumference, the area of the circle is seen 
to have been of direct interest to the 
authors of the \suls. 
There is no indication in the \suls\ of  their being aware of  
the ratio of
the circumference to the diameter being the same as the ratio of the area
to the square of the radius; no occasion seems to have presented itself 
that would inspire a comparison of the two ratios. The issue of area, which was involved
in the construction of altars,  is treated independently. 
There were {\it citi}s (fire platforms) constructed in the
shape of a chariot wheel, a circular trough etc. with stipulated areas, which 
motivated the issue of how to transform 
a square into a circle having the same area.

\medskip
\noindent {\bf Transforming a square into a circle}
 
\medskip
\bsl\ describes a  procedure to produce a circle with the same area as a given
square, which  goes as follows: take a string with length half the diagonal of the square, 
and stretch from the
centre across a side of the square, viz. PS as in Figure~1,  and draw the
circle including a third of the 
extra part stretching outside the square, viz. PR as in the figure with QR=$\frac 13$ PS, as
the radius.

\begin{figure}[ht!]
\begin{center}
\includegraphics[width=65mm]{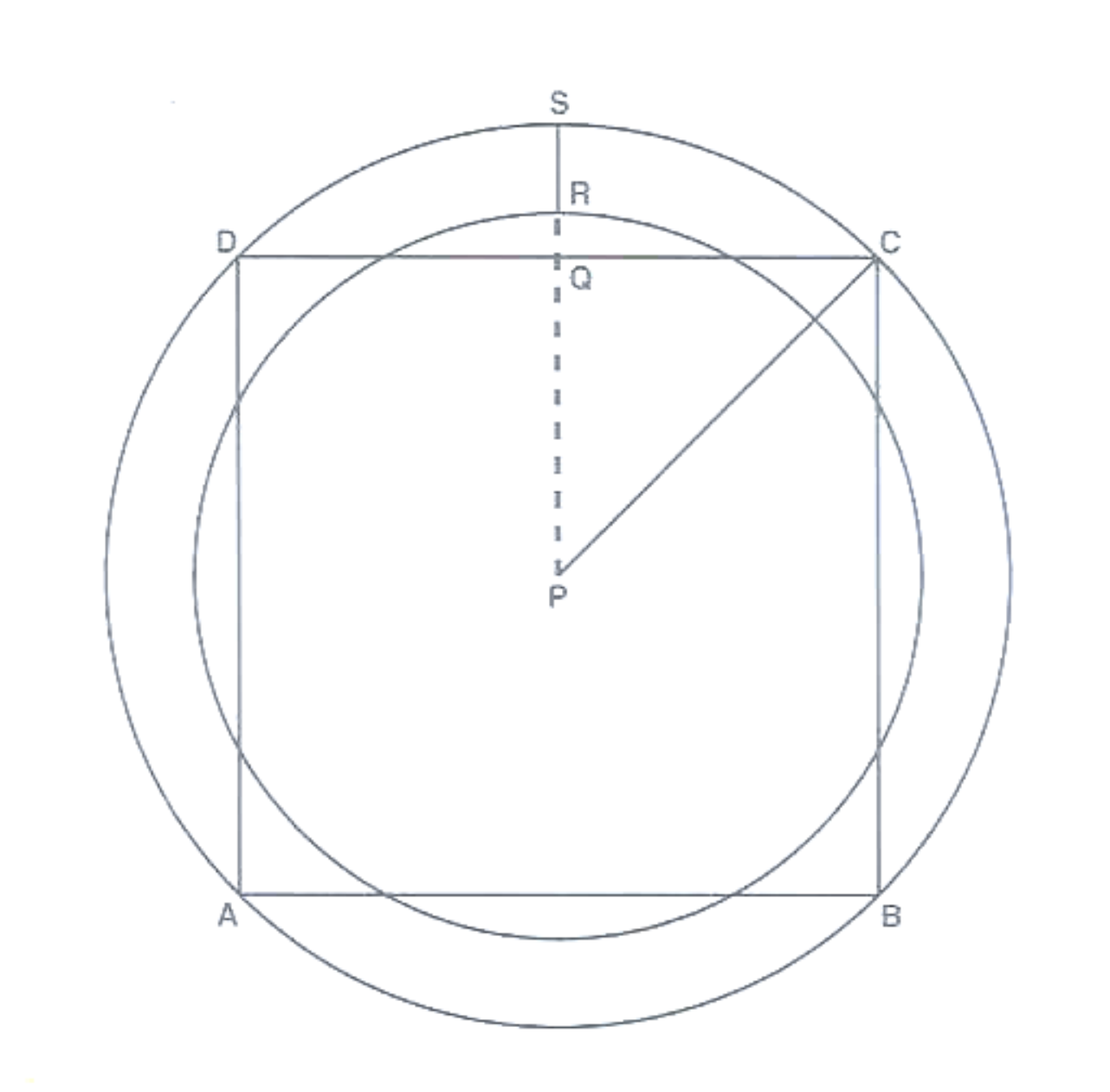}\\
{Circling the square, \bsl\ \sul}
\end{center}
\end{figure}

For a square with side $2a$ this radius works out to 
be $a+\frac 13 (\sqrt 2 -1)a=(2+\sqrt 2)a/3 $. 
For the unit square the area of this circle works out to  $1.01725...$, about
$1.7 \%$ more than the correct value 1.
If one is to compute the value of $\pi$ with $(2+\sqrt
2)a/3 $ as the radius of the circle corresponding to a square with side  
$2a$, it works out to be $3.0883\dots $, in place of
$3.1415\dots$. It 
should be borne in mind that what they had was a procedure for
producing the circle and not a numerical value for $\pi$; the latter
had not emerged as a  concept, and they were not trying to compute such
a ratio. The comparison, here and in similar contexts below, is
only for facilitation in overall comprehension of the relative
values.

\medskip
The \asl\ \sul\ gives the same construction for the circle, as \bsl. 
\msl\ \sul\ is seen to provide another construction for the circle with the area of
the given square. The following interpretation of a verse in \msl\ \sul\
was introduced in \cite{D} by this author. For convenience I shall discuss the 
verse and background around it, after first describing the procedure, according
to the interpretation.  The steps involved 
are illustrated by Figure~2.

\begin{figure}[ht!]
\centering
\includegraphics[width=65mm]{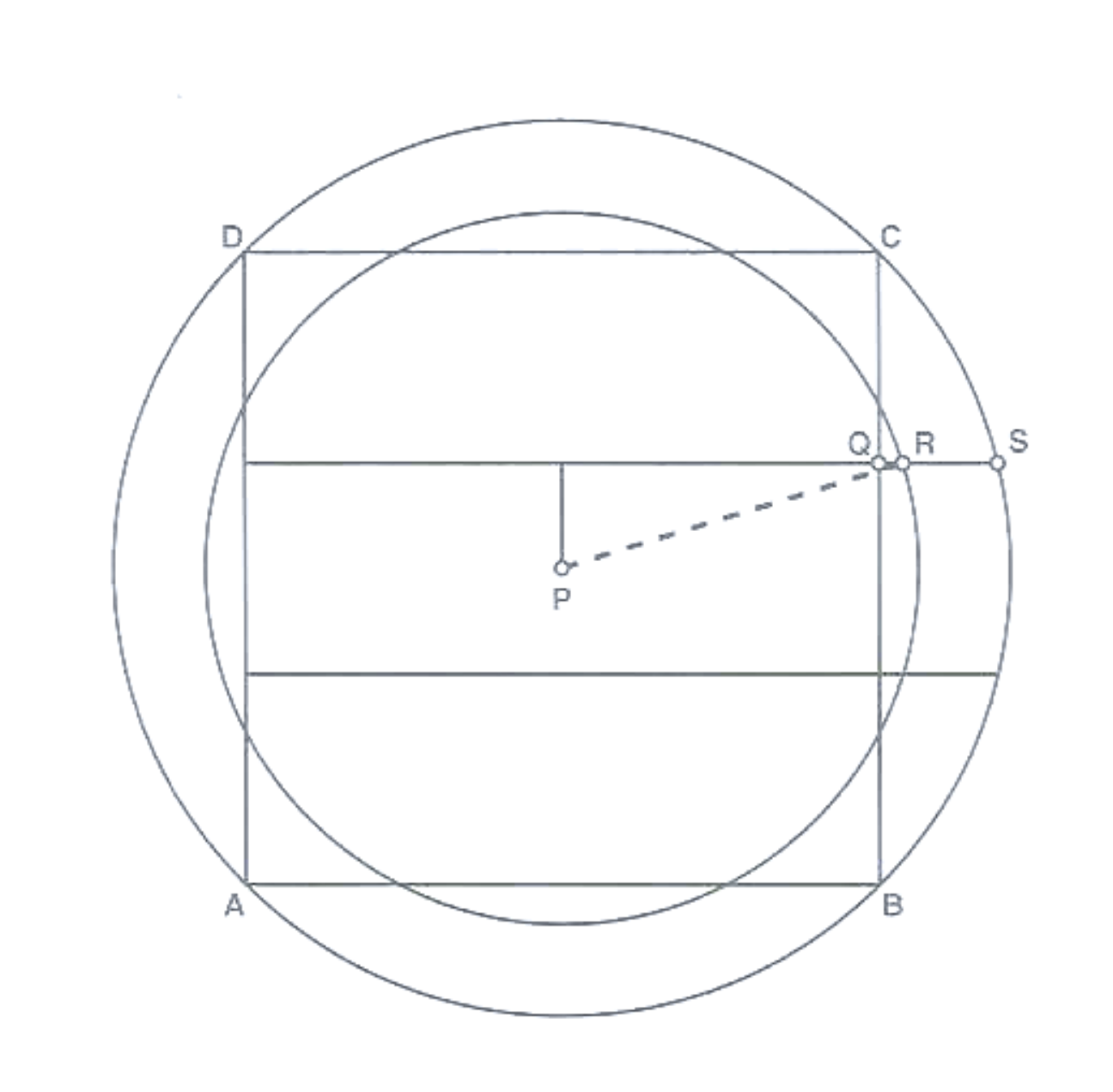}\\
{Circling the square, \msl\ \sul}
\end{figure}

\begin{quote}
Draw the lines dividing the square into 3 equal
strips. Produce one of these lines to meet the circle passing through
the vertices of the square. On the segment of the line that is outside
the square and inside the circle, viz. QS as in  Figure~2, take the
point at a distance $\frac 15$th of the length of the segment, from
the square, viz. the point R in the figure, with QR=$\frac 15$QS. The
circle with PR as the radius, where P is the centre of the square, is
given as the desired circle, with area equal to that of the original
square.     
\end{quote}

For a square of
side length $2$ the length of the segment 
between the square and the outer circle is seen to be $\displaystyle \left(\frac 
{\sqrt{17}}3 -1\right )$, so the radius $r$ of the circle is given by 
$$r^2=\left \{1+\frac 15 \left (\frac {\sqrt {17}}3 -1\right )\right \}^2+\frac 19.$$ 
For the unit square the area of this circle works out to $0.9946\dots$, 
a substantially more accurate value compared to the earlier one, the
error involved being just about $\frac 12 \%$ (which is now on the other side). 
The value of $\pi$ in this case works out to  $3.1583\dots$.

The verse in question, from  \msl\  \sul\  is:

\begin{center}
{\tt catura\'sra\d m navadh\=a kury\=addhanu\d hko\d tyastridh\=atridh\a $|$\\

utsedh\=atpa$\tilde {\tt n}$cama\d m lumpetpur\={\i}\d se\d neha t\=avatsama\d m $||$}

\smallskip
\hfill (\msl\ \sul\ 11.15)
\end{center}

There seems to be considerable confusion about, and discomfort with,
this verse in literature.  In \cite{SB} it is suggested that possibly ``squares are drawn without any 
mathematical significance".  In \cite{Ku} there is an interpretation, the conclusion of which manifestly wrong. There is another interpretation in \cite{RCG-m}, concluding the value
of $\pi$ according to the s\u tra to be $\frac {25}8$, but it may be seen that in the interpretation 
the first line of the verse plays no role at all, while that of the second is quite ad hoc. 
Appropriate transliteration seems to be at the heart of the issue. 

I propose the following translation: 

\begin{quote}
Divide the square into nine parts, (by) dividing the sides into three (equal) parts each.  
Mark a fifth of the part jutting out (of the square) and cover (the corresponding circle with centre at the origin) with loose earth. 
\end{quote}

The meaning of this (according to my interpretation) is described above, with the help of  Figure 2. It would be out of place to go into the linguistic details with regard to the interpretation. I shall instead focus on highlighting two reasons for which the present interpretation ought to be appropriate. Firstly, as noted above, it leads to a significantly improved result; this could not be a mere coincidence.   Secondly if one is to try to read their mind on how they might have attempted to remedy a perceived discrepancy in a known result, the construction seems to arise as a natural development: In the first place, it is reasonable to suppose in this respect that over a period it had been realised that the circle produced by the \bsl\ construction was slightly larger than it should be. Since taking a point on the bisector of the square along a side did not seem to work, they chose to consider trisectors of square. So far there is no divergence in various interpretations. The crucial, and distinctive point in the above interpretation is that they picked a point on the trisecting line, which is actually natural in the context of the comparison with the \bsl\ construction, but does not seem to have been taken note of by the earlier authors. Furthermore in analogy with the earlier construction a point had to be picked up on the segment of the trisector jutting out from the square. In the earlier construction one third of the jutting out part was added to get the radius of the desired circle, and it may be noted here that though the circle through the vertices of the square finds no mention in the verse, it would be lurking in the their minds, in the context of the \bsl\ construction. The fraction $\frac 15$th of the extra part is then likely to have been based on an ad hoc observation that the resulting line segment for the radius is slightly smaller than the \bsl\ construction, as was desired.

Clearly, the \msl\ construction is the result of keen attempts to
improve upon the original \bsl\ construction, through refinement of
the overall scheme. How
the specific details  were conceived and how, and to what extent, it was
confirmed to be more accurate, remains unclear.

It may be mentioned here that there were also other constructions adopted, in the broader Vedic community; while indeed the Vedic civilization shared a certain common body of knowledge, there are many variations in the individual \suls\ adopted by different sub-communities. A lesser known \sul\ by the name {\it Maitr\a ya\d n\={\i}ya},
which is akin to \msl\ \sul\ (in that it belongs to the same {\it sa\. mhita}), 
gives a construction for circling the square which involves
taking the radius of the desired circle to be $\frac 9{16}$ times the
side of the square; see~\cite{RCG}. For a unit square the area of the
resulting circle 
turns out to be  $0.9940\dots$, comparable to the one above, with
 $3.1604\dots $ as the corresponding value of $\pi$. It may
be recalled here that the Egyptians took the area of a circle of
radius $r$ to be $(\frac {16r}9)^2$, to  which the above value, for the
reverse process, corresponds exactly. 

\medskip
\noindent {\bf Squaring the circle}

\medskip
The converse problem of  ``squaring the circle'', viz. finding a square
with the same area as a given circle\footnote{This should not be confused with 
the problem in Greek geometry of finding a square with {\it precisely} the 
same area through a {\it ruler and compass} construction. The context of 
the s\u trak\a ras was entirely different, and their objective would have been only to find a square with the 
area of the circle, within the levels of accuracy they were used to, or desired. 
They may have liked to find a geometric construction, like 
by ruler and compass, but that was not the thrust. The problem in their 
perspective involved finding such a square by whatever available means.}  
is also considered in the \suls.
\bsl\  gives the following expression for the ratio of the side of
the square to the diameter of the circle (the original description is
in words): 
 $$~~~~~~~~~~~ {\frac 78}+\frac 1{8\times 29}- \frac 1{8\times 29
   \times 6}+\frac 1{8\times 29\times 6 \times 8}. ~~~~~~~~~~~(1)$$ 
For the circle
with unit 
radius, the area according to this works out to $3.0883...$, a little more than
$98.3\%$ of the actual value. 
It may be noted that in this case also the error is about 1.7 percent, in the
opposite direction. 
It could  not be a coincidence (as has been noted also by earlier authors - see [19], [17]), that the errors in the two prescriptions, corresponding to mutually opposite operations, while substantial, are quite matching in their order and opposite in the orientation.  It suggests that for want of a geometrical procedure in the reverse direction (unlike for transforming a square into a circle) they obtained it through inversion of the previous ratio, in some way which is not entirely clear so far (see below). To get the inverse of $(2 +\sqrt 2)/3$ they sought out a value of  $\sqrt 2$,  as a familiar fraction.

\medskip
\noindent{\bf The square root of 2}

\medskip
Three of the four \suls, 
\bsl, \asl ~and \ksl, give the following expression for $\sqrt 2$ (in
words):
$$~~~~~~~~~~~ 1+\frac 13 + \frac 1{3\times 4}-\frac 1{3\times 4 \times
  34}.~~~~~~~~~~~ (2)$$
In decimal expansion the value of the expression is $1.4142157\dots$. 
This is remarkably close to the actual value $1.4142136\dots$, and this fact 
has been a subject of much laudatory comment in literature. It 
may be recalled  in this context that Babylonians also had  a value, about a
thousand years earlier, describing $\sqrt 2$ in the sexagesimal system,
which works out to $1.4142129...$ (see \cite{FR}, for instance).  Various aspects
including the presentation of the number as above and the substantial relative difference
of the values (including the side of the error), rule out any organic link between the
values. There has been considerable speculation and discussion on how the \sul\
value of $\sqrt 2$ may have been arrived at; we shall however not digress to these
details here (see for instance \cite{D} and other references cited there).  

As noted above the motivation for finding a value for $\sqrt 2$ would have come  from the problem of computing the inverse of $(2+\sqrt 2)/3$, viz. for
obtaining the formula~(1).  This numerical value of $\sqrt 2$ is not involved
elsewhere in the \suls.  
In other contexts they are seen to use only the geometric
form of $\sqrt 2$ as the diagonal of the unit square, which in fact
went by the special name {\it dvikara\nd \={\i}}. 

How the 
inversion would have been effected, using the value of $\sqrt 2$ as
above, has been discussed by Thibaut~\cite{T} and also other later
authors. The older explanations, however, presuppose considerable dexterity 
on the part of the s\u trak\a ras in dealing with fractions, for which there 
is no corroborative evidence, and are thus unsatisfactory. In a recent 
paper Kichenassamy~\cite{Ki} has proposed a resolution of the issue which is more 
in tune with the \suls\ ethos; the paper also 
discusses at length the inadequacies of the earlier arguments. 

It would also be worthwhile to note here another \sul\ construction 
which relates in a way to properties of the circle. \bsl\ \sul\ describes a construction of a square 
which involves drawing a perpendicular to a given line, say L, at a point P on L, by
drawing circles with centres at points on either side of P on L at equal distance, with
radii larger than the distance, and joining the points of intersection of the two 
circles (in very much the same way as taught in schools today). Underlying 
the construction is the realisation that the line joining the two points of 
intersection of two circles meets the line joining their centres orthogonally; though
the construction involves the principle for circles of equal radius, it seems reasonable
to assume that they were aware of this ``orthogonality principle" in that generality. In 
most constructions requiring perpendiculars, they were however produced using 
the converse of the Pythagoras theorem\footnote{The theorem named after Pythagoras
has been known in India at least since the time of the earliest \bsl\ \sul\ (ca. 800 BCE), where an
explicit statement of it is found. The converse of the theorem, namely that a triangle in which 
the square of one of the sides equals the sum of squares of the other two sides, was used 
extensively for producing perpendiculars (see \cite{D} for details).}, rather than the construction as above 
(implemented in a certain way, the former turns out to be simpler than the latter; 
see~\cite{D}). 

Let me conclude this section on the \suls\ with the following comment. There has
been a tendency with regard to \suls\ to assume that the s\u trak\a ras lay great 
store on accuracy. While the value of $\sqrt 2$ does seem like an example 
of this, a careful reading of the \suls\ shows that  high degree of accuracy was not seen as a 
primary objective. In many
contexts, alternate values or constructions are described, which are of a crude variety,
alongside some relatively accurate ones, which shows that in their overall conception,
the benefits meant to accrue from the ritual performances would not be seriously 
affected if approximate procedures were adopted. Where accuracy was pursued,
it seems to be the result of keen academic enquiry, rather than an imperative arising
from practical issues of the time, or the 
philosophical framework involved. On the other hand a supposition as above does them a 
disservice in the context of the less accurate values such as in circling a square. 
Mathematics of ancient cultures needs to be understood and appreciated  in their specific context, 
and not judged through generalised abstract tests. The issue of circling the square arose
for instance from the desire of having a fire platform with the same area, that would
not bring with it an intrinsic demand for high degree of accuracy, and it is incorrect to wonder why their
value of the area is not accurate - {\it it was not meant to be}. 

\bigskip
\noindent{\bf {\large The Jaina tradition} }

\medskip
Apart from the Vedic religion (if one may call it that) Jainism and Buddhism
flourished during the first millennium BCE (and later during certain periods). 
There was a long tradition among the {Jaina}s of engagement with
mathematics, as is evidenced from various compositions that have come 
down to us. As for Buddhism, though certain constructions involved in 
Buddhist pursuit,  called {\it Mandala}s 
involve intricate designs which seem mathematically significant, no textual 
composition involving mathematical concepts has come down to us. 

The motivation of the Jainas for mathematics did not came, per se,  from any rituals, which they
indeed abhorred, but from contemplation of the cosmos, of which they
had evolved an elaborate and unique conception of their own. In the {\it 
  Jaina} cosmography the world is supposed to be an infinite flat plane, with
concentric annular regions surrounding an innermost circular region
with a diameter of 100000 {\it yojana}s\footnote{{\it yojana} was a measure of 
length, of the order of 15 to 20 kilometres, with local variations.}, known as the {\it Jambudv\={\i}pa}
(island of Jambu, that corresponded to the Earth), and the annular 
regions alternately consist of water and land, and the width
of each successive ring being twice that of the previous one. 
The geometry of the circle played an 
important role in the overall discourse, even when the scholars
engaged in it were more of philosophers than
mathematicians. Unfortunately, many historical and chronological
details of the {Jaina} 
tradition are uncertain (even more so than the Hindu tradition) and 
have been a subject of speculation. Definitive references to
the properties of the circle known in the {Jaina} tradition can be
found in the work in fourth or fifth century (\cite{P}, page 59). It
has however been mentioned by Datta in~\cite{D-j} that they are found
in {\it  Tattv\a rth\a  
  dhigama-s\u tra-bh\a \S ya}, a philosophical work of Um\a sv\a ti, who
is supposed to have lived around 150 BCE according to the {\it \'Svet\a
  mbara} tradition and in the second century CE according to the {\it
  Digambara}  tradition\footnote{{\it \'Svet\a
  mbara} and  {\it
  Digambara} are two ancient branches of Jainism, with certain differences 
  both in terms of their philosophy as also practices in everyday life.}. Datta~\cite{D-j} suggests that   Um\a sv\a ti was
probably not the discoverer of the formulae and they would have been
known centuries before him, and discusses some evidence in this
respect. Saraswati Amma (\cite{SA}, page 63) attributes the basic 
formulae to   
{\it S\u ryapraj\~n\a pti}, a composition  which is believed to be from the fifth
century~BCE. 

In the {Jaina} tradition the departure from old belief of 3 as the
ratio of the circumference to the diameter is quite pronounced; 
{\it S\u ryapraj\~n\a pti} records the then traditional value 3 for
it, and discards it in 
favour of $\sqrt {10}$. The {Jaina}s also knew the ratio of the
area of the 
circle to the square of the radius to be the same number as the ratio of the circumference
to the diameter. In fact
they had the formula directly relating the area to the circumference, 
that the former is a fourth of the product of the circumference and
the diameter of the circle, which in particular readily implies the equality
of the two ratios as above.  Incidentally, $\sqrt
{10}$ which is about $3.16227\dots$,  
may be seen to be a better approximation for $\pi$ compared
to the \bsl\ construction, 
involving an error of only about  $\frac 23$~per cent. 

 The value was very
convenient to the Jaina theologian mathematicians, in their computations. For example, in {\it Jambudv\={\i}pa praj$\tilde n$apti} the value
the circumference of {\it Jambudv\={\i}pa}, a circle of diameter $100,000$ {\it yojana}, is computed, with $\sqrt {10}$ as the value for the ratio of the circumference to the diameter, by computing 
the square-root of $10^{11}$. 

This value for $\pi$ 
was used for over a thousand years, even after better
values were known; indeed so routine was its use in the {Jaina} 
texts that it is often known as the {Jaina} value for $\pi$. 
The value was also adopted in 
{\it Pa$\tilde n$casiddh\a ntika}, in the Siddh\a nta tradition sometime 
 during 1st to 6th century, and by Brahmagupta in the 7th century. 

\medskip
There has been some speculation on the origin of $\sqrt {10}$ as a
value for $\pi$. One 
explanation, attributed to Hunrath, goes as follows (\cite{SA}, page 65): 
The square of the side of a regular 12-sided polygon inscribed in a circle of unit 
radius is seen to be 
$(1 - \frac {\sqrt 3}2)^2+(\frac 12)^2$ and with the
choice of $\frac 53$ as an approximation for $\sqrt 3$ it works out to  $\frac {\sqrt {10}}{12}$;
thus the perimeter of a regular 12-gon is about $\sqrt {10}$ times its principal diagonals, but 
a little less.  This may have
 inspired  the formula for the circumference of the circle. The explanation 
however involves many assumptions about which there is little evidence. An interested
reader may also consult  \cite{RCG-j} for various other explanations. 

V\={\i}rasena, a Jaina mathematician from the 8th century
states:  
\begin{center}
{\tt vy\a sam \sd o\dd a\s agu\nd ita\md\ \sd o\dd a\s a sahita\md\ trir\u par\u
pairbhakta\d m $|$

vy\a sam trigu\nd ita\md\ s\u k\sd m\a dapi tadbhavet s\u k\sd ma\d m $|$}

\end{center}

\hfill (\d Sa\d tkha\nd \d d\a gama, Vol. IV. page 42)

\noindent A routine  translation of this would be as follows (slightly simplified from~\cite{S}): 

\begin{quote}
Sixteen times the diameter, together with 16, divided by
113 and thrice the diameter is a very fine value (of the circumference). 

\end{quote}

There is something strange about the formula (with the interpretation as above), that it prescribes
``together with~16'' - surely it was known to the author that the
circumference is proportional to the diameter and that adding 16,
independent of the diameter,  
would not be consistent with this. It seems reasonable however to suppose 
that the author meant 
$3+\frac {16}{113}=\frac
{355}{113}$ to be the factor by which to multiply the diameter to get the circumference, which is indeed a good approximation, as the author
stresses with the phrase ``s\u kshmadapi  s\u
kshamam'' (finest of the fine !).\footnote{There does not seem to be much of 
scope for attributing the issue to corruption in the course of transmission at some
level; there are however issues of grammar and interpretation involved, and this accepted 
translation may be flawed; it is possible, for instance, that {\it \sd o\dd a\s a sahita\md} ("together with sixteen"), 
which is the culprit, has the role of emphasising  that while dividing by 
$113$, one is to divide the previous product, which involved $16$ - thus 
``together with 16" is not about adding 16, but reiterating that the following 
division by $113$ is to be subjected to the output together with the earlier 16. 
While this indeed does seem odd in what is intended to be a formula, it need not be ruled out, given that the 
direct interpretation is odd anyway, and the author obviously would not have meant it. In part, the somewhat curious presentation may have been been the result of needs of versification.}  In China this approximation for $\pi$ was given by
Chong-Zhi (429-500). Its value is 3.1415929... in place of
3.1415926...,  
accurate to 6 decimals. 

In {\it Trilokas\a ra}, which is another account of the {Jaina} 
scholarship, composed by Nemicandra, who lived around 980 CE, also one
finds another value for the ratio $\pi$, apart from $\sqrt {10}$: it
is  the value $(\frac {16}9)^2$, that we saw from the 
{\it Maitr\a ya\n \i ya} \sul\ (shared also with the
Egyptians). This may suggest a relation with the Hindu
tradition, but the time gap is rather intriguing. 

Apart from the circle as a whole, the {Jaina} mathematicians were also interested
in the interrelation between the arcs of a circle and the corresponding chords. This is related to their conception of the 
geography of {\it Jambudv\={\i}pa}, including various regions,
mountains etc.  Um\a sv\a ti notes various relations
between the length $c$ of a chord, the height  $h$  of the corresponding ``arrow''
(viz. the segment joining the midpoint of the chord to the midpoint of
the arc)  and the diameter $d$ 
of the circle. One of the relations noted is $$c=\sqrt{4h(d-h)};$$ various 
other forms which are equivalent to this one algebraically, from a modern 
point of view, are also presented. An interested reader may consult \cite{S}, \cite{D-j} and 
also \cite{RCG-bow} 
for further discussion on this issue; the last two references have some details 
also of analogous formulae from other ancient cultures.  

There is also an interesting 
formula for the length of the minor arc (the smaller of the
two arcs cut out by the chord), say $a$, as 
$$a=\sqrt{6h^2+c^2}$$ (with notation as above). As can be seen, such 
a relation does not actually hold exactly. It may be noted that in the
special case when
the chord is a diameter, so that the arc is a semicircle, the equation
corresponds to the ratio of the length of the semicircle to the radius
being $\sqrt {10}$; so the relation holds with {\it their} value for
$\pi$. As we go to small arcs however the assertion goes quite off the
mark. Surprisingly however, the formula continued to be part of {
  Jaina} literature all the way, including the famous mathematical 
work {\it Ga\nd ita s\a ra sa\. ngraha} of Mah\a v\={\i}ra in 850 (Ch. VII, verse 73$\frac 12$; see 
\cite{P}, page 469). The formula
also appears  in {\it Trilokas\a ra} of Nemicandra (see \cite{D-n}),
who was mentioned above. 

Given a chord of a circle, apart from the length of the arc segment
one may also ask about the area cut out by the chord (with the minor
arc).   {\it Ga\nd ita s\a ra sangraha} gives the value of the area to
be  $\frac
14 \sqrt {10}ch$, where $c$ is the length of the chord and $h$ is the
height over the chord (length of the arrow). The formula is also found
in {\it Trilokas\a ra} of Nemicandra. This formula also holds strictly only 
for a semicircular segment with $\pi$ in place of $\sqrt {10}$ but
diverges from the actual value for smaller arcs. 


\medskip
A different formula for the area of the segment cut out by the chord
is given in {\it Tri\s atika} of  \'Sridhara (ca. 750)\footnote{Though until some time 
there had been an argument over when he lived and his background, there is now a general consensus that 
\'Sridhara is from the 8th century, and was Jaina, at least during the time of his writings - his 
mathematical work is seen to be consistent with this.}, and also quoted in some later 
works, including Bh\a skara (see below for more about him). 
It may be worth mentioning here 
that \'Sridhara  does not quite seem to fit 
in the astronomer mathematicians tradition - his known works deal
exclusively with mathematics, and he is well-known for his procedure for
solving a quadratic equation.  According to \'Sridhara 
the area $A$ of 
the segment between a chord and the corresponding arc is given by  
$$A=\frac {\sqrt {10}}3 {\displaystyle\left (\frac{h(c+h)}2\right )};$$
Clearly  $\sqrt {10}$ here is meant to be for the ratio $\pi$. 

It would seem
that many formulae for arc segments cut out by chords were written
down by extrapolating relations that 
were noted for the case of the semicircle to a general arc segment; if
they had some (heuristic) reasoning for it, it is not found 
recorded. 
From a historical point of view this highlights the
difficulties faced by the ancient mathematicians in grasping the lengths
of arcs and the areas bounded by them, and their endeavour to get around
the difficulties, before the ideas of trigonometry, and then calculus  
emerged.

\bigskip
\noindent{\bf {\large {\= A}ryabha\d ta and the astronomical tradition}}

\medskip
{{\= A}ryabha\d ta}, born in 476 CE (as has been indicated by the author in his work {\it \=Aryabha\d t\={\i}ya}), was the pioneer     
of what is termed as the {\it siddh\a nta} tradition, of astronomer mathematicians in India that flourished for almost
eight hundred years, until Bh\a skar\a ch\a rya in the 12th century, and even beyond, 
and in turn led to the Kerala school of 
mathematics.  While the
tradition has some manifest linkages with the older Hellenistic mathematical  
astronomy, after the early influences it seems to have charted a course of its own. 
Many new  
mathematical ideas were developed, both in response to the theoretical
demands in the study of astronomy, and also in pursuit of pure
mathematical thought. In particular a deeper 
understanding of the circle evolved,  both in terms of geometry and
trigonometry. In his work {\= A}ryabha\d t\={\i}ya we find the following:

\begin{center}
{\tt Caturadhikam \s atama\sd tagu\nd am dv\a \sd a\s tistath\a\ sahasr\a \nd am
$|$

ayutadvayavi\sd kambhasy\a sanno v\r ttapari\nd ahah $||$}

\smallskip
\hfill (Ga\nd itap\a da 10, in {\it \A ryabha\d t\={\i}ya})
\end{center}

\begin{quote}
The circumference of a circle with diameter twenty thousand is
approximately a hundred and four times eight, and sixty-two thousand [viz. 62832].

\end{quote}
 
\noindent This gives the value of $\pi$ as approximately $3.1416$, which indeed coincides 
with the correct value of $\pi$ truncated at $4$~decimal places. 
It may be recalled that in Greek astronomy,  Ptolemy
had the value, in sexagesimal expression, which corresponds to $3.14166...$. 
There is no direct information on how {\= A}ryabhata arrived at the value. One
may anticipate that, like in similar instances in other cultures, the value was 
obtained through repeated application of the
formula $$S_{2n}=\sqrt{\frac {S_n^2}4+\left (1-\sqrt{\frac {4r^2-S_n^2}4}\right )^2},$$
where $S_n$ is the side of the regular $n$-gon inscribed in a
circle of unit radius. The formula follows from the ``Pythagoras
theorem'', which, as noted earlier, had been known in India since the \bsl\ \sul\ (8th century BCE) and is also stated in {\=Aryabhat\={\i}ya} (in Ga\d nitap\=ada - 17). It is 
 suggested by Gane\s a, a sixteenth century commentator of {\=
 A}ryabha\d t\={\i}ya, that an inscribed  polygon with 384 sides  
was used as an approximation for the circle, and the
above formula was used, starting with a hexagon (for which the side 
coincides with the radius of the circle), until reaching the
polygon with the number of sides $384=6\times 2^6$. The choice of 
20,000 as the measure for diameter is readily seen to facilitate computation
of square-roots in integral values; the values would have been rounded, up or down, to 
integer values at various stages of application of the above formula, and the square 
root computed using the well-known procedure for the purpose, that is attributed to \A ryabha\d ta.
It may be noted that the value of $\pi$ as above is slightly greater than the actual 
value, despite its representing the perimeter of an inscribed regular polygon, due to 
rounding up at some stages. 

In {\= A}ryabhatiya we also have the trigonometric sine functions.\footnote{While the Greeks did trigonometry with chords, it was in India
that the trigonometry in terms of the half-chords originated.} 
{{\= A}ryabha\d t\={\i}ya} (499 CE) provides a sine table, in a verse, for 
angles upto $90^\circ$ that are multiples of $3^\circ 45'$ ($24$th part of 
the right angle): taking  the circle whose circumference
is 21600 (equal to the total measure of the circumference in minutes), the differences between the values of half-chords corresponding to angles that are successive multiples of $3^\circ 45'$
are recounted sequentially; the radius of the circle, which features as the total of the differences recounted, is 3438. A similar table also appears in {\it Panca-siddh\a ntik\a}
an older composition from the early centuries of CE in which the value of the radius 
involved is 120. Once such tables were available, the
lengths of circular arcs could be calculated using the sine table (for the specific 
values), without 
recourse to any special formula as in {Jaina} mathematics. There were also interpolation
methods for dealing with intermediate angles. 

Apart from the sine tables there was also a curious approximate formula
for the sine function in vogue in the Siddh\a nta tradition. It is generally 
attributed to Bhaskara I (7 th century CE), being part of his {\it Mah\a bh\a skar\={\i}ya}, 
but is also found independently in the contemporaneous work {\it Br\a hmaspu\d ta siddh\a nta} of 
Brahmagupta. In the modern notation the formula may be stated as 
$$\sin \theta = \frac {4\theta (180-\theta)}{40500-\theta (180-\theta)},$$
where $\theta $ is the angle measured in degrees. The formula is seen to be remarkably
accurate, involving an error of less than $1\%$, except for very small angles. It is unclear
how such a formula was derived. (see \cite{RCG-sine} and \cite{vB} for further details in this respect). 

Knowledge of various properties of the circle and trigonometry gradually 
became crucial part of learning in the Siddh\a nta tradition, being a 
prerequisite for pursuing mathematical astronomy. The tradition sustained itself, 
though perhaps somewhat feebly during certain periods than others, and
individual exponents made fresh contributions to knowledge, apart from 
carrying forward the body of knowledge that was getting built. We shall
not go into the finer historical details in this respect here. 
Bh\a skara II, from the 12th century, (also known as Bh\a skar\a c\a rya, Bh\a skara the teacher) 
is considered  the last major exponent from the tradition. Apart from mastering
the knowledge flowing in the tradition, Bh\a skara made substantial contributions of his own
in various respects. 
By his time, the
attendant mathematics, especially arithmetic and geometry, that went with 
mathematical astronomy, had acquired a wider appeal, and applicability,  in the society.  
Bh\a skar\a c\a rya composed a comprehensive work, {\it Siddh\a nta \'Siroma\d ni} which, in the 
tradition of Siddh\a nta works, had a chapter devoted the mathematical topics 
as above, called {\it L\={\i}lavat\={\i}}. The latter however acquired a life of its own, and 
a reputation as a mathematical work, with
large number of copies being produced. It served as a textbook of 
mathematics for several centuries, in  a large part of India. Specifically with 
regard to the circle I will only recall the following  (approximate) formula from {\it L\={\i}lavati} for the length of an arc of a circle; the formula itself may be seen to be related to Bh\a skara I's formula
for the sine function, when expressed in radians:
$$a=\frac p2 -\sqrt{\frac {p^2}4-\frac {5p^2c}{4(c+4d)}}=\frac p2 \left (1-\sqrt {1-\frac{5c}{c+4d}}\right ),$$
where $p$ denotes the circumference (perimeter) of the circle, and the other
notation is as above, namely $a$ is the length of the arc, $c$ is the length of the chord,
and $d$ is the diameter of the circle; the first expression as above is akin to the way it is given in the
original verse and the second is a simplification. 
 A more integral view is seen to have
evolved with regard to geometry of the circle and trigonometry.

\bigskip
\noindent{{\large \bf The Kerala School}}

\medskip
We conclude this article with a few observations on the Kerala school
in the context of the above theme. The school originated with the work of 
M\a dhava in the second half of the 14th century, and flourished, as a teacher-student
continuity, with multiple names involved during some periods, for about 250 years. 
They took remarkable strides towards calculus, introducing techniques involving
infinitesimals, and in particular had obtained Gregory-Leibnitz series for the arctan 
function and the Newton series for the sine function (over two centuries before
their European counterparts). 
We shall not go into a detailed
discussion on the mathematics from the Kerala school, which has been a
subject of much study in recent years. The interested reader is
referred to \cite{GJ}, \cite{P}, and \cite{RS}.  

Determining accurate values for $\pi$, which is something that concerns our 
theme here, seems to have been a passion for the school. In particular 
the following remarkably close approximation to $\pi$ is credited to
M\a dhava, by \'Sankara V\=ariar (1556), in {\it Kriy\=akramakar\={\i}} (cf. \cite{P}): 
the measure of the circumference in a circle of diameter 
$900,000,000,000$ is $2,827,433,388,233.$
Thus $$\pi =\displaystyle \frac
{2,827,433,388,233}{900,000,000,000}=3.141592653592\dots ,$$
in place of  $3.1415926535 89\dots$, accurate to 11 decimals, when rounded. 
As  the series expansion

\medskip
\centerline{Circumference = $ 4\,\mbox {\rm diameter}\  ( 1-\frac 13 + \frac 15
  + \dots )$}  
  
  \medskip
 \noindent that they had obtained converges very slowly, 
and hence  not useful in getting good approximations for $\pi$. 
Madhava had introduced an ingenious device to get over this difficulty, 
called {\it antya samsk\a ra}, ``the end correction''. 
With $S_n$ as the sum of the series truncated at the $n$th term,
he introduced sequences $a_n$ such that the sequence $S_n+a_n$
converges faster. 
The third, the final one that was recorded, produces the
sequence 
$$S_n+ (-1)^{n-1}\displaystyle \frac{n^2+1}{4n^3+5n}.$$
The 50 th term of this is accurate to 11 decimals. 

\medskip
\noindent {\bf Note:} The papers of R.C. Gupta cited here are also available in 
the compilation of {\it Ga\d nit\a nanda}, edited by K. Ramasubramanian, Published by the Indian Society for History of Mathematics (ISHM), 2015.

\vskip8mm
\begin{flushleft}
Department of Mathematics\\
Indian Institute of Technology\\
Powai, Mumbai, 400005 

\smallskip
E-mail: {\tt sdani@math.iitb.ac.in}
\end{flushleft}

\end{document}